\documentclass[preprint,oneside,letterpaper]{revtex4}
\usepackage{hyperref}
\usepackage{amsmath,amsfonts}
\usepackage{graphicx}
\usepackage{subfigure}
\usepackage{enumerate}
\usepackage{psfrag}
\usepackage{tikz}
\newtheorem{lemma}{Lemma}
\newtheorem{thm}[lemma]{Theorem}

\newcommand{\qed}{$\hfill \square$}
\begin{document}

\title{Local matching indicators for transport with concave costs}

\author{Julie Delon}
\email{julie.delon@enst.fr}
\affiliation{LTCI CNRS, T{\'e}l{\'e}com ParisTech}
\author{Julien Salomon}
\email{salomon@ceremade.dauphine.fr}
\affiliation{Universit{\'e} Paris IX/CEREMADE}
\author{Andre{\u\i} Sobolevski{\u\i}}
\email{ansobol@mccme.ru}
\address{A. A. Kharkevich Institute for Information Transmission Problems, Moscow, Russia}
%\address{UMI 2615 CNRS ``Laboratoire J.-V.~Poncelet''}
\thanks{This work is supported by ANR
  through grant ANR-07-BLAN-0235 OTARIE
  (\url{http://www.mccme.ru/~ansobol/otarie/}); AS thanks the Ministry
  of National Education of France for supporting his visit to the
  Observatoire de la C{\^o}te d'Azur
  and to Russian Fund for Basic Research for the partial
  support via grant RFBR 07--01--92217-CNRSL-a.}

\begin{abstract}
   In this note, we introduce a class of indicators that enable to
  compute efficiently optimal transport plans associated to arbitrary
  distributions of $N$ demands and $N$ supplies in $\mathbf{R}$ in the
  case where the cost function is concave. The computational
  cost of these indicators is small and independent of $N$. A
  hierarchical use of them enables to obtain an efficient algorithm.
\end{abstract}

\maketitle
\section{Introduction}
It is well known that transport problems on the line involving convex cost
functions have explicit solutions, consisting in a monotone
rearrangement. Recently, an efficient method has been introduced to
tackle this issue on the circle~\cite{DSS}. In this note we introduce an algorithm that 
enables to tackle optimal transport problems on the line (but actually
also on the circle) with
concave costs. Our
algorithm complements the method suggested by McCann~\cite{McCann}.  McCann
considers general real values of supply and demand and shows how the
problem can be reduced to convex optimization somewhat similar to the
simplex method in linear programming.  Our approach as presented here
is developed for the case of unit masses and is
closer to the purely combinatorial approach of~\cite{Aggarwal}, but extends it to
a general concave cost function. The extension to integer masses will be presented in~\cite{dss}.\\
The method we propose is based on a class of local indicators,
that allow to detect consecutive points that are matched in an optimal
transport plan. Thanks to the low number of evaluations of the cost
function required to apply the indicators, we derive an algorithm that
finds an optimal transport plan in $n^2$ operations in the worst
case. In practice, the computational cost of this method appears to behave linearly with respect
to $n$. \\
Since the indicators apply locally, the algorithm can be massively
parallelized and also allows to
treat optimal transport problems on the circle. In this way, it extends 
the work of Aggarwal {\it et al.}~\cite{Aggarwal} in which cost functions
have a linear dependence in the distance.
\section{Setting of the problem}
For $N_0\in \mathbf{N}^*$, consider $P=(p_i)_{i=1,...,N_0}$ and  $Q=(q_i)_{i=1,...,N_0}$ two sets of
points in $\mathbf{R}$ that represent respectively demand and supply locations.
The problem we consider in this note consists in minimizing the transport cost
\begin{equation}\label{eq:4}
C(\sigma) = \sum_{i,j} c(p_i, q_{\sigma(i)}),
\end{equation}
where $\sigma$ is a permutation of $\{1,...,N_0\}$. This permutation
forms a \emph{transport plan}.\\ 
We focus on the case where the function $c$ involves a concave
function as stated in the next definition.
\begin{def}\label{def:cost}
  The \emph{cost function} in~(\ref{eq:4}) is defined on $\mathbf{R}$ by
$ c(p, q) = g(|p - q|)$ with $p, q \in \mathbf{R}$, where $g(\cdot)$ is a concave non-decreasing
  real-valued function of a real positive variable such that $g(0) :=
  \lim_{x\to 0} g(x) \ge -\infty$. %
\end{def}
Some examples of such costs are given by $g(x) =
\log (x)$ with $g(0) = -\infty$, and $g(x) = \sqrt x$ or $g(x) = |x|$ is with
$g(0) = 0$.\\
Finally, we denote by $\sigma^\star$ the permutation associated to a given
optimal transport plan between $P$ and $Q$: for all permutation $\sigma$ of
$\{1,...,N_0\}$, $C(\sigma^\star)\leq C(\sigma).$

\section{Chains}\label{sec:noncross}
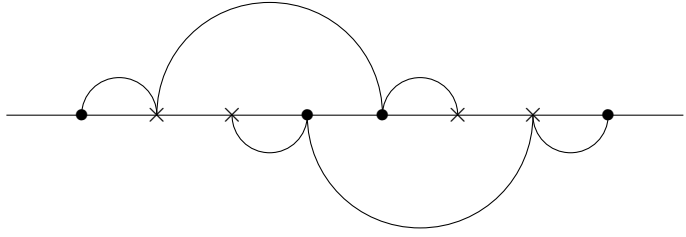
\begin{figure}
\begin{center}
\begin{tikzpicture}
\draw (0,0) --   (9,0) ;

\draw (1,0) node {$\bullet$} ;
\draw (4,0) node {$\bullet$} ;
\draw (5,0) node {$\bullet$} ;
\draw (8,0) node {$\bullet$} ;

\draw (2,0) node {$\times$} ;
\draw (3,0) node {$\times$} ;
\draw (6,0) node {$\times$} ;
\draw (7,0) node {$\times$} ;

\draw (2,0) arc  (0:180: .5) ;
\draw (5,0) arc  (0:180:1.5) ;
\draw (6,0) arc  (0:180: .5) ;

\draw (4,0) arc  (0:-180: .5) ;
\draw (7,0) arc  (0:-180:1.5) ;
\draw (8,0) arc  (0:-180: .5) ;
\end{tikzpicture}
\caption{Example of chains.}\label{exchaine}
\end{center}
\end{figure}

In this section, we present a way to build a particular partition of
the set $P\cup Q$.\\
%\subsection{Non-crossing rule}
Consider two pairs of matched points $(p_i,q_{\sigma^\star(i)})$ and 
$(p_{i'},q_{\sigma^\star(i')})$, say e.g. $p_i\leq
q_{\sigma^\star(i)}$, $p_{i'}\leq q_{\sigma^\star(i')}$. It is easy to
prove that the  
following alternative holds:
\begin{enumerate}
\item $[p_i,q_{\sigma^\star(i)}]\cap [p_{i'},q_{\sigma^\star(i')})]=
  \emptyset$,
\item $[p_i,q_{\sigma^\star(i)}]\subset
  [p_{i'},q_{\sigma^\star(i')})]$ or
  $[p_{i'},q_{\sigma^\star(i')})]\subset[p_i,q_{\sigma^\star(i)}]$.
\end{enumerate}
This remark is a direct consequence of the concavity of the cost
function and is often denominated as "the non-crossing rule"
\cite{Aggarwal,McCann}. In the next section, we show
how it allows decompose the
initial situation in sub-problems where supply and demand points are
alternated. \\
%\subsection{Chains}
Because of the non-crossing rule in an optimal plan there are as many
supply points as demand points between any pair of matched points $p_i$
and $q_{\sigma(i)}$.  For a given demand point $p_i$, define its {\it left
neighbor} $q'_i$ as the nearest supply point on the left of $p_i$ such that
the numbers of supply and demand points between $q'_i$ and $p_i$ are equal;
define the {\it right neighbor} $q''_i$ of $p_i$ in a similar way.  Then define
a {\it chain} as a maximal alternating sequence of supply and demand
points ($p_{i_1}, q_{j_1}, p_{i_1}, ..., q_{j_k}$) or ($q_{j_1}, p_{i_2},
..., p_{i_{k + 1}}$) such that each $q_{i_l}$ is the right neighbor of
$p_{i_l}$ and the left neighbor of $p_{i_{l + 1}}$.  An extension of the
proof of Lemma 3 of \cite{Aggarwal} shows that the collection of chains forms a
partition of the set $P \cup Q$. 
An simple example of such a partition is shown on~Fig.~\ref{exchaine}.
 Note that construction of this
collection only depends on relative positions of supply and demand
points and does not involve any evaluation of the cost function.  It
can be done in $O(N_0)$ operations.\\
The non-crossing rule implies that all matched pairs of points in an
optimal transport plan must belong to the same chain.  We therefore
restrict ourselves in the sequel, without loss of generality, to the
case of a single chain 
\begin{equation}\label{alter}
p_1 < q_1 < ... < p_i < q_i < p_{i + 1} < q_{i
+ 1} < ... < p_N < q_N, 
\end{equation}
for $N\in\mathbf{N}^*$ and keep these last notations throughout the rest of this paper.
\section{Local matching indicators}\label{sec:LMI}
Thanks to the non-crossing rule, one knows that there exists at least two consecutive
points $(p_i,q_i)$ or $(q_i,p_{i+1})$ that are matched in any optimal
transport plan. Starting from this remark, we take advantage of the
structure of a chain to introduce
a class of indicators that enable to detect a priori such pairs of points.
\begin{def}\label{def:LMI}%(Local Matching Indicators of order $k$)
We define 
$$
I^p_k(i)=c(p_i    ,q_{i+k})+\sum_{\ell=0}^{k-1}c(p_{i+\ell+1},q_{i+\ell})-\sum_{\ell=0}^{k}c(p_{i+\ell},q_{i+\ell}),
$$
where $k,i$ are such that $1\leq k\leq N-1$ and $1\leq i\leq N-k$, and
$$
I^q_k(i)=c(p_{i+k+1},q_{i})+\sum_{\ell=1}^{k}c(p_{i+\ell},q_{i+\ell})-\sum_{\ell=0}^{k}c(p_{i+\ell+1}  ,q_{i+\ell}),
$$
for $k,i\in \mathbf{N}$, such that $1\leq k\leq N-2$ and $1\leq i\leq N-k-1$.
\end{def}
The interest of these functions lies in the next result.

\begin{figure}
\begin{center}
\begin{tikzpicture}
\draw (0,0) --   (7,0) ;

\draw (1,0) node {$\bullet$} ;
\draw (3,0) node {$\bullet$} ;
\draw (5,0) node {$\bullet$} ;

\draw (2,0) node {$\times$} ;
\draw (4,0) node {$\times$} ;
\draw (6,0) node {$\times$} ;

\draw (4,0) arc  (0:180: .5) ;
\draw (5,0) arc  (0:180:1.5) ;

\draw (0+8,0) node {\Huge $\leq$} ;

\draw (0+9,0) --   (7+9,0) ;

\draw (1+9,0) node {$\bullet$} ;
\draw (3+9,0) node {$\bullet$} ;
\draw (5+9,0) node {$\bullet$} ;

\draw (2+9,0) node {$\times$} ;
\draw (4+9,0) node {$\times$} ;
\draw (6+9,0) node {$\times$} ;

\draw (3+9,0) arc  (0:180: .5) ;
\draw (5+9,0) arc  (0:180: .5) ;

\draw (8,-2) node {\Huge $\Downarrow$} ;

\draw (0+4.5,-4) --   (7+4.5,-4) ;

\draw (1+4.5,-4) node {$\bullet$} ;
\draw (3+4.5,-4) node {$\bullet$} ;
\draw (5+4.5,-4) node {$\bullet$} ;

\draw (2+4.5,-4) node {$\times$} ;
\draw (4+4.5,-4) node {$\times$} ;
\draw (6+4.5,-4) node {$\times$} ;

\draw (4+4.5,-4) arc  (0:180: .5) ;
\end{tikzpicture}
\caption{Schematic representation of the result stated in Thm.~\ref{ThLMI} in
  the case where $k=1$.}\label{pic:ThLMI}
\end{center}
\end{figure}
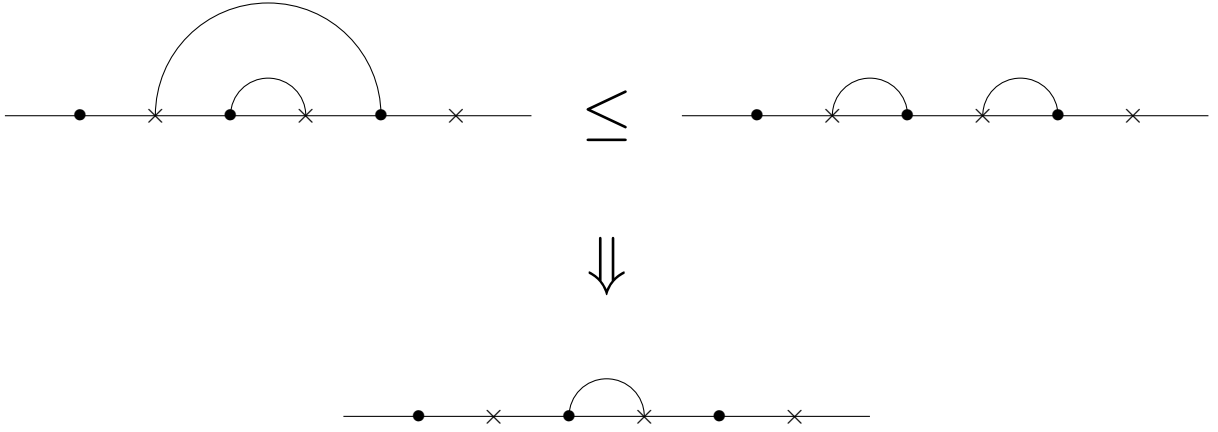

\begin{thm}\label{ThLMI}%[Negative Local Matching Indicators of order $k$]
Let $k_0\in \mathbf{N}$ with $1\leq k_0\leq N-1$ and $i_0\in \mathbf{N}$ (resp. $i_0'\in \mathbf{N}$), such
that $1\leq i_0\leq N-k_0$ (resp. $1\leq i_0'\leq N-k_0-1$).\\  
 Assume that 
\begin{enumerate}
\item  $I^p_k(i) \geq 0$ for $k=1,...,k_0-1$, $1\leq i \leq N-k$,\label{hyp1}
\item  $I^q_k(i')\geq 0$ for $k=1,...,k_0-1$, $1\leq i'\leq N-k-1$, \label{hyp2}
\item  $I^p_{k_0}(i_0) < 0$ (resp. $I^q_{k_0}(i_0')< 0$). \label{hyp3}
\end{enumerate}
Then any permutation $\sigma^\star$ associated to an optimal transport
plan satisfies
$\sigma^\star(i)=i-1$ for $i=i_0+1,...,i_0+k_0$
(resp. $\sigma^\star(i)=i$ for $i=i_0+1,...,i_0+k_0$).
\end{thm}
In practice,
these indicators allow to find pairs of neighbors that are
matched in an optimal transport plan.
This result is illustrated on~Fig.~\ref{pic:ThLMI}.
\\
Before giving the proof, we state a basic result.
\begin{lemma}\label{lemm:basic}
We keep the previous notations. Define 
$$
\varphi_{k,i}^p(x,y)=g(x+y+p_{i+k}-q_i)+\sum_{\ell=0}^{k-1}c(p_{i+\ell+1},q_{i+\ell})-g(x)-g(y)-\sum_{\ell=1}^{k-1}c(p_{i+\ell},q_{i+\ell}), 
$$
for $k,i\in \mathbf{N}$, such that $1\leq k\leq N-1$ and $1\leq i\leq N-k$, and
$$
\varphi_{k,i}^q(x,y)=g(x+y+p_{i+k+1}-q_{i})+\sum_{\ell=1}^{k}c(p_{i+\ell},q_{i+\ell})-g(x)-g(y)-\sum_{\ell=1}^{k-1}c(p_{i+\ell+1}  ,q_{i+\ell}),
$$
for $k,i\in \mathbf{N}$, such that $1\leq k\leq N-2$ and $1\leq i\leq N-k-1$.
Both functions $
\varphi_{k,i}^p(x,y)$ and $
\varphi_{k,i}^q(x,y)$ are decreasing with respect to each of their two variables.
\end{lemma}
This lemma is a direct consequence of the concavity of the function
$g$.
We are now in the position to give the sketch of the proof of Theorem \ref{ThLMI}.\\
{\bf Proof of Theorem~\ref{ThLMI}: }
We consider the case where $I^p_{k_0}(i_0)< 0$. The case
$I^q_{k_0}(i_0')<0$ can be treated the same way.\\
The proof consists in proving that Assumptions (\ref{hyp1}--\ref{hyp3}) imply that
neither demand nor supply points located between $p_{i_0}$ and
$p_{i_0+k_0+1}$ can be matched with points located outside this interval,
i.e. that the set $\mathcal{S}_{i_0}=\{ p_i, i_0+1\leq i \leq i_0+k_0 \}\cup \{q_i,
i_0\leq i \leq i_0+k_0-1\}$ is stable by an optimal transport plan. In
this case, the result follows from Assumption (\ref{hyp1}--\ref{hyp2}).\\
Suppose that $\mathcal{S}_{i_0}$ is not preserved by an optimal
transport plan $\sigma^\star$. Three cases can occur:
\begin{enumerate}[a)]
\item\label{c1} There exists $i_1\in \mathbf{N}$, such that $1\leq i_1\leq i_0$ and
  $ i_0 \leq \sigma^\star(i_1)\leq i_0+k_0 -1$ and there exists $i_1'\in
  \mathbf{N}$, such that $\sigma^\star(i_1)+1 \leq i_1'\leq i_0+k_0$ and
  $ i_0+k_0 \leq \sigma^\star(i_1')\leq N$.
\item\label{c2} There exists $i_2\in \mathbf{N}$, with $i_0+1\leq i_2\leq i_0+k_0$ such that
  $ 1 \leq \sigma^\star(i_2)\leq i_0-1$.
\item\label{c3} There exists $i_2\in
  \mathbf{N}$, with $i_0+k_0 < i_2 \leq  N$ such that
  $ i_0 \leq \sigma^\star(i_2)< i_0+k_0$.
\end{enumerate}
We first prove that Case \ref{c1}) cannot occur.\\
In Case \ref{c1}), one can assume without loss of generality that $\sigma^\star (i_1)$ is the largest
index such that $1\leq i_1\leq i_0$, $ i_0 \leq \sigma^\star(i_1)\leq
i_0+k_0 -1$ and that $i_1'$ is the smallest index such that
$\sigma^\star(i_1)+1 \leq i_1'\leq i_0+k_0$, $ i_0+k_0 \leq \sigma^\star(i_1')\leq N$ 
. With such assumptions, the (possibly empty) subset $\{ p_i, \sigma^\star (i_1)+1 \leq i \leq i_1'-1 \}\cup \{q_i,
\sigma^\star (i_1)+1\leq i \leq i_1'-1\}$ is stable by
$\sigma^\star$. Because of Assumptions (\ref{hyp1}--\ref{hyp2}), no 
nesting can occur in this subset, and $\sigma^\star(i)=i$ for
$i=\sigma^\star(i_1)+1,...,i_1'-1$.\\
On the other hand, since $\sigma^\star$ is supposed to 
be optimal, one has:
$$
c(p_{i_1},q_{\sigma^\star(i_1)})+c(p_{i_1'},q_{\sigma^\star(i_1')})+\sum_{i=\sigma^\star(i_1)+1}^{i_1'-1}c(p_i,q_i)
\leq
c(p_{i_1},q_{\sigma^\star(i_1')})+\sum_{i=\sigma^\star(i_1)}^{i_1'-1}c(p_{i+1},q_i).
$$
Thanks to Lemma \ref{lemm:basic},
one deduces from this last inequality that:
$$
c(p_{i_0},q_{\sigma^\star(i_1)})+c(p_{i_1'},q_{i_0+k_0})+\sum_{i=\sigma^\star(i_1)+1}^{i_1'-1}c(p_i,q_i)\\
\leq
c(p_{i_0},q_{i_0+k_0})+\sum_{i=\sigma^\star(i_1)}^{i_1'-1}c(p_{i+1},q_i),
$$
and then:
\begin{eqnarray}
c(p_{i_0},q_{\sigma^\star(i_1)})+\sum_{i=i_0}^{\sigma^\star(i_1)-1}c(p_{i+1},q_i)
+c(p_{i_1'},q_{i_0+k_0})+\sum_{i=i_1'}^{i_0+k_0-1}c(p_{i+1},q_i)
\nonumber\\
+\sum_{i=\sigma^\star(i_1)+1}^{i_1'-1}c(p_i,q_i)\leq
c(p_{i_0},q_{i_0+k_0})+\sum_{i=i_0}^{i_0+k_0-1}c(p_{i+1},q_{i}).\ \ \label{ineq:1}
\end{eqnarray}
According to Assumption (\ref{hyp1}),
$I^p_{\sigma^\star(i_1)-i_0}(i_0)\geq 0$ and
$I^p_{i_0+k_0-i_1'}(i_1')\geq 0$, so that:

$$\sum_{i=i_0}^{\sigma^\star(i_1)}c(p_i,q_i)\leq
c(p_{i_0},q_{\sigma^\star(i_1)})+\sum_{i=i_0}^{\sigma^\star(i_1)-1}c(p_{i+1},q_{i})$$
$$\sum_{i=i_1'}^{i_0+k_0} c(p_i,q_i)\leq
c(p_{i_1'},q_{i_0+k_0})+\sum_{i=i_1'}^{i_0+k_0-1}c(p_{i+1},q_i).
$$
Combining these last inequalities with (\ref{ineq:1}) one finds that:
$$\sum_{i=i_0}^{i_0+k_0}c(p_i,q_i) \leq c(p_{i_0},q_{i_0+k_0})+\sum_{i=i_0}^{i_0+k_0-1}c(p_{i+1},q_i),
$$
which contradicts Assumption (\ref{hyp3}). \\
Let us now prove that Cases \ref{c2}) and
\ref{c3}) contradict the assumptions.
Cases \ref{c2}) and \ref{c3}) can be treated in the same
way. Consider Case \ref{c2}). Without loss of generality, one can
assume that $i_2$ is the smallest index such that $i_0+1\leq i_2 \leq
i_0+k_0$ and $\sigma^\star(i_2)\leq i_0-1$. Because there are
necessarily as much demands as supplies between $q_{i_0}$ and
$p_{i_2}$, there exists one and only one index $i_2'$ such that
$i_0\leq \sigma^\star(i_2')\leq i_2-1$ and $1\leq i_2'\leq
i_0$. Consequently, the (possibly empty) subsets $\{ p_i, i_0 +1 \leq i \leq \sigma^\star(i_2')  \}\cup \{q_i,
i_0   \leq i \leq \sigma^\star(i_2')-1\}$ and $\{ p_i,
\sigma^\star(i_2')+1 \leq i \leq i_2-1 \}\cup \{q_i,
\sigma^\star(i_2')+1\leq i \leq i_2 -1\}$ are stable by an optimal
transport plan. Because of Assumptions (\ref{hyp1}--\ref{hyp2}), no
nesting can occur in these subsets, and $\sigma^\star(i)=i-1$ for
$i=i_0+1,...,\sigma^\star(i_2')$ and $\sigma^\star(i)=i$ for
$i=\sigma^\star(i_2')+1,...,i_2-1$. \\
On the other hand, since $\sigma^\star$ is supposed
to be optimal, one has
\begin{eqnarray*}
c(p_{i_2},q_{\sigma^\star(i_2)})+c(p_{i_2'},q_{\sigma^\star(i_2')})+\sum_{i=i_0+1}^{\sigma^\star(i_2')}c(p_{i},q_{i-1})
+\sum_{i=\sigma^\star(i_2')+1}^{i_2-1}c(p_{i},q_{i})
\\\leq c(p_{i_2'},q_{\sigma^\star(i_2)}) + \sum_{i=i_0+1}^{i_2}c(p_{i},q_{i-1}).
\end{eqnarray*}
Thanks to Lemma~\ref{lemm:basic},
one deduces from this last inequality that:
\begin{eqnarray}
c(p_{i_2},q_{\sigma^\star(i_2)})+c(p_{i_0},q_{\sigma^\star(i_2')})+\sum_{i=i_0+1}^{\sigma^\star(i_2')}c(p_{i},q_{i-1})
+\sum_{i=\sigma^\star(i_2')+1}^{i_2-1}c(p_{i},q_{i})\nonumber
\\\leq c(p_{i_0},q_{\sigma^\star(i_2)}) + \sum_{i=i_0+1}^{i_2}c(p_{i},q_{i-1}).\label{ineq3}
\end{eqnarray}
Because the cost is supposed to be increasing with respect to the
distance, one finds that
$c(p_{i_0},q_{\sigma^\star(i_2)})\leq c(p_{i_2},q_{\sigma(i_2)})$, so that (\ref{ineq3}) implies:
%\begin{eqnarray*}
$$c(p_{i_0},q_{\sigma^\star(i_2')})+\sum_{i=i_0+1}^{\sigma^\star(i_2')}c(p_{i},q_{i-1})
+\sum_{i=\sigma^\star(i_2')+1}^{i_2-1}c(p_{i},q_{i})
\leq \sum_{i=i_0+1}^{i_2}c(p_{i},q_{i-1}),$$
%\end{eqnarray*}
and then:
\begin{eqnarray}
c(p_{i_0},q_{\sigma^\star(i_2')})+\sum_{i=i_0+1}^{\sigma^\star(i_2')}c(p_{i},q_{i-1})
+\sum_{i=\sigma^\star(i_2')+1}^{i_2-1}c(p_{i},q_{i})+\sum_{i=i_2+1}^{i_0+k_0}c(p_{i},q_{i-1})\nonumber\\
\leq \sum_{i=i_0+1}^{i_0+k_0}c(p_{i},q_{i-1}).\label{ineq4}
\end{eqnarray}
According to Assumption (\ref{hyp1})
$I^p_{\sigma^\star(i_2')-i_0}(i_0)\geq 0$, so that:

$$\sum_{i=i_0}^{\sigma^\star(i_2')}c(p_i,q_i)\leq
c(p_{i_0},q_{\sigma^\star(i_2')})+\sum_{i=i_0}^{\sigma^\star(i_2')-1}c(p_{i+1},q_i).$$
Combining these last inequalities with (\ref{ineq4}) one finds that:
$$\sum_{i=i_0}^{i_0+k_0}c(p_i,q_i) \leq c(p_{i_0},q_{i_0+k_0})+\sum_{i=i_0}^{i_0+k_0-1}c(p_{i+1},q_i),
$$
which contradicts Assumption (\ref{hyp3}).\\
We have then shown that 
neither demand nor supply points located between $p_{i_0}$ and
$q_{i_0+k_0+1}$ can be matched with located outside this interval. The
set $\mathcal{S}_{i_0}$ is then stable by an optimal transport
plan. According to Assumption (\ref{hyp1}--\ref{hyp2}), no nesting can
occur in $\mathcal{S}_{i_0}$. The result follows.  \qed
\section{An algorithm for balanced chains}\label{sec:algo}
The recursive use of our indicators is on the basis of the next algorithm. \\
{\bf Algorithm:}
Set $\mathcal{P}=\{p_1,...,p_N,q_1,...,q_N\}$, $\ell^p=(1,...,N)$, $\ell^q=(1,...,N)$, and $k=1$.\\
While $\mathcal{P}\neq \emptyset $ and $k\leq N-1$ do
\begin{enumerate}
 \item Compute $I^p_k(i)$ and $I^q_k(i')$ for $i=1,...,N-k$ and $i'=1,...,N-k-1$.\label{costlystep}
 \item Define $$\mathcal{I}_k^p=\{i_0, 1\leq i_0\leq N-k ,I^p_k(i_0)<0 
  \},$$ $$\mathcal{I}_k^q=\{i_0, 1\leq i_0\leq N-k-1,I^q_k(i_0)<0 \},$$ 
  and do \label{noncostlystep}
 \begin{enumerate}
  \item If $\mathcal{I}_k^p=\emptyset$ and $\mathcal{I}_k^q= \emptyset$, set $k=k+1$.
  \item Else do \label{removingstep}
  \begin{itemize}
   \item for all $i_0$ in $\mathcal{I}_k^p$ and for $i=i_0  +1,...,i_0+k$, do
   \begin{itemize}
    \item define $\sigma^\star(\ell^p_i)=\ell^q_{i-1}$,
    \item remove $\{p_{\ell^p_i},q_{\ell^q_{i-1}}\}$ from $\mathcal{P}$,
    \item remove $\ell^p_i$ and $\ell^q_i$ from $\ell^p$ and $\ell^q$ respectively.
   \end{itemize}
   \item for all $i_0'$ in $\mathcal{I}_k^q$ and for $i=i_0'+1,...,i_0'+k$, do
   \begin{itemize}
    \item define $\sigma^\star(\ell^p_i)=\ell^q_i$,
    \item remove $\{p_{i},q_{i}\}$ from $\mathcal{P}$,
    \item remove $\ell^p_i$ and $\ell^q_i$ from $\ell^p$ and $\ell^q$ respectively.
   \end{itemize}
   \item set $N=\frac 12 Card(\mathcal{P})$, and rename the points in
     $\mathcal{P}$ such that
     $\mathcal{P}=\{p_1,...,p_N,q_1,...,q_N\}$,  
$$p_1<q_1<...<p_{i}<q_{i}<p_{i+1}<q_{i+1}<...<p_{N}<q_{N}. $$
   \item set $k=1$.
  \end{itemize}
 \end{enumerate}
\end{enumerate}
If $k=N-1$, for $i=1,...,N$ set $\sigma^\star(\ell^p_i)=\ell^q_i$.\\

\begin{centering}
\begin{figure}[h!]
\psfrag{g}[c][t]{In-line evaluations of $g$}
\psfrag{N}[c][b]{Number of pairs of points $N$}
%\psfrag{data1}[l][c]{: Best case $N-1$}
\psfrag{data1}[l][c]{\tiny: $g(x)=|x|^{10^{-3}}$, $\alpha=1.18$}
\psfrag{data2}[l][c]{\tiny: $g(x)=\sqrt{|x|}$, $\alpha=1.87$}
\psfrag{data3}[l][c]{\tiny: $g(x)=|x|^{1-10^{-3}}$, $\alpha=2$}
\psfrag{data4}[l][c]{\tiny: Worst case $(N-1)^2$, $\alpha=2$}
\psfrag{100}[l][c]{\small $100$}
\psfrag{150}[c][c]{\small $150$}
\psfrag{200}[c][c]{\small $200$}
\psfrag{250}[c][c]{\small $250$}
\psfrag{300}[c][c]{\small $300$}
\psfrag{350}[c][c]{\small $350$}
\psfrag{400}[c][c]{\small $400$}
\psfrag{450}[c][c]{\small $450$}
\psfrag{500}[c][c]{\small $500$}
\psfrag{10}[c][c]{\small $10$}
%\psfrag{alpha}[l][c]{$\alpha$}
\includegraphics[width=.8\textwidth,height=.5\textwidth]{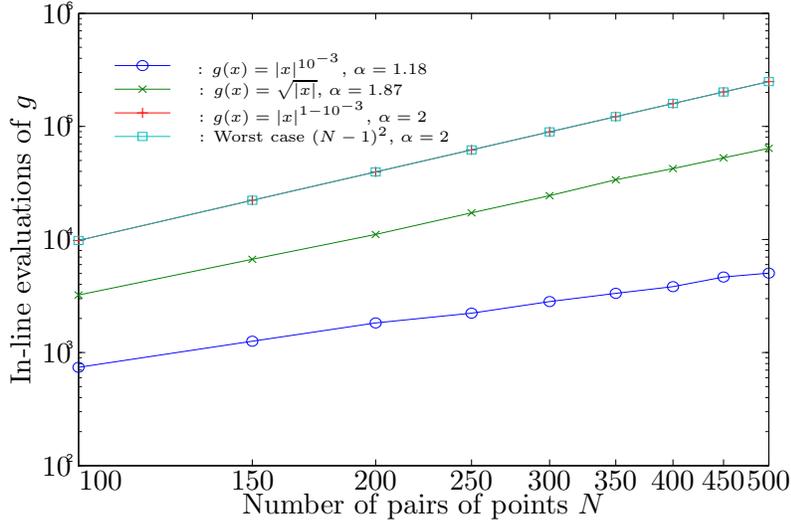}
\caption{Number of in-line evaluations with respect to the number of
  pairs, for various cost functions. The number $\alpha$ is the slope of the log-log graphs.}\label{Ec}
\end{figure}
\end{centering}

To test the efficiency of our algorithm, we have applied it to an
increasing number $N$ of pairs of points. 
For a fixed value of $N$,
100 samples of points have been chosen
randomly in $[0,1]$, and the mean of the number of evaluations of $g$ has been
computed. The results are shown on~Fig.~\ref{Ec}. \\
The best case consists in finding a negative indicator at
each step, and the worst corresponds to the case where all the
indicators are positive. These two cases require respectively $N-1$
and $(N-1)^2$ evaluations of $g$.\\

\end{document}